\documentclass{article}
\usepackage[utf8]{inputenc}
\usepackage{mathtools}

\usepackage{amssymb}
\usepackage{amsthm}
\usepackage[left=1in, right=1in, top=1in, bottom=1in]{geometry}
\usepackage{verbatim}
\normalfont\upshape
\theoremstyle{plain}
\usepackage{fancyhdr}
\usepackage{setspace}
\usepackage{hyperref}

\newtheorem*{theorem*}{Theorem}
\newtheorem{theorem}{Theorem}
\newtheorem*{lemma*}{Lemma}
\newtheorem{lemma}{Lemma}
\newtheorem*{corollary*}{Corollary}

\newtheorem{proposition}{Proposition}
\theoremstyle{remark}
\newtheorem*{remark*}{Remark}

\newtheorem*{fact*}{Fact}
\newtheorem*{proposition*}{Proposition}
\newtheorem*{definition*}{Definition}

\newtheorem*{remarks}{Remarks}

\title{A Strong Splitting of the Frobenius Morphism on the Algebra of Distributions of $SL_2$}
\author{\small Gus Lonergan}
\date{}

\begin{document}

\maketitle

\begin{abstract}\noindent Let $p$ be a prime number, and let $Dist(SL_2)$ be the algebra of distributions, supported at $1$, on the algebraic group $SL_2$ over $\mathbb{F}_p$. The Frobenius map $Fr:SL_2\to SL_2$ induces a map $Fr:Dist(SL_2)\to Dist(SL_2)$ which is in particular a surjective algebra homomorphism. In this note, we construct a (unital) section of this map, whenever $p\geq 3$. The main ingredient of this construction is a certain congruence modulo $p^3$, reminiscent of the congruence $\binom{np}{p}\equiv n\mod p^3$.\end{abstract}

\section{Introduction}

Fix a prime number $p$. Let $G$ be an affine algebraic group defined over $\mathbb{F}_p$. Following e.g. \cite{J}, we consider its \emph{algebra of distributions} $H=Dist(G)$. This is an augmented Hopf algebra analogous to the universal enveloping algebra in characteristic $0$. The analogy is strong when $G$ is simply-connected and semisimple, in the sense that its category of finite-dimensional representations is equivalent to the category $Rep(G)$ of finite-dimensional (algebraic) representations of $G$. However, structurally it is very different from the universal enveloping algebra, being (for instance) not finitely generated. Indeed, it may in fact be regarded as a `divided power' version of the universal enveloping algebra.\

Nonetheless, its structure may with some effort be studied. One main tool is the Frobenius morphism $Fr: H\to H$ (induced by the usual Frobenius endomorphism of $G$). $Fr$ is a surjective (augmented Hopf algebra) endomorphism of $H$, whose kernel is equal to the augmentation of a certain finite-dimensional augmented subalgebra $H_1$ of $H$. In fact, $H_1$ is nothing more than the algebra of distributions of the kernel of the Frobenius morphism on $G$. Taking distribution algebras of kernels of higher and higher powers of Frobenius, we get an exhaustive filtration:\begin{align*}\mathbb{F}_p=H_0\subset H_1\subset H_2\subset\ldots \end{align*}of $H$. Each $H_i$ is an augmented subalgebra of $H$, of dimension $p^{i.dim(G)}$; and the Frobenius endomorphism induces surjections $H_i\to H_{i-1}$ (for $i\geq1$), and the identity map $H_0\to H_0$.\

In \cite{G1}, \cite{G2}, the authors construct (for $G=SL_2$ in \cite{G1}, and for any $G$ simply-connected semisimple in \cite{G2}) a certain non-unital splitting of $Fr:H\to H$. This is a non-unital map of algebras (but neither augmented nor Hopf) $\phi:H\to H$ such that $Fr\circ \phi=Id$. In their splitting, the image of $1$ is equal to a certain idempotent element of $Dist(T)$ (for a choice $T$ of $\mathbb{F}_p$-split maximal torus of $G$), whose effect on any finite-dimensional representation $V$ of $G$ is to project to the sum $V_0$ of all $T$-weight subspaces of weights divisible by $p$. Consequently, their splitting amounts to giving $V_0$ the structure of representation of $G$. In other words, the construct a functor $Rep(G)\to Rep(G)$, lying over the functor $Rep(T)\to Rep(T)$ given by $V\to V_0$.\

The aim of this paper is, in the case $G=SL_2$, to upgrade this to a \emph{unital} splitting $\theta=\theta_T$. Equivalently, for any $r\in\mathbb{F}_p$, consider the functor $Rep(T)\to Rep(T)$ which sends the finite-dimensional representation $V$ of $T$ to the sum $V_r$ of its $T$-weight subspaces of weights congruent to $r\mod p$; then we lift this to a functor $Rep(G)\to Rep(G)$. We achieve this by giving explicit generators and relations for $H$ which make it rather clear. Namely, for a $\mathbb{F}_p$-split maximal torus $T$ we choose a standard basis $\{e,f,h\}$ of $Lie(G)=\mathfrak{sl}_2(\mathbb{F}_p)$ such that $h$ spans the Lie algebra of $T$, and we have:

\begin{theorem}\label{thm:1}$H$ is generated by the elements $e,e^{(p)},e^{(p^2)},\ldots$ and $f,f^{(p)},f^{(p^2)},\ldots$ subject to the relations:

\begin{enumerate}\item $[X_k,e^{(p^k)}]=2e^{(p^k)}, [X_k,f^{(p^k)}]=-2f^{(p^k)}$,\
\item $[X_k,e^{(p^{k+n})}]=0=[X_k,f^{(p^{k+n})}]$,\
\item $[e^{(p^k)},e^{(p^{k+n})}]=0=[f^{(p^k)},f^{(p^{k+n})}]$,\
\item $[e^{(p^k)},f^{(p^{k+n})}]=(-1)^{n}(f^{(p^k)})^{p-1}(f^{(p^{k+1})})^{p-1}\ldots (f^{(p^{k+n-1})})^{p-1}(X_k+1)$,\\
 $[e^{(p^{k+n})},f^{(p^{k})}]=(-1)^{n}(X_k+1)(e^{(p^k)})^{p-1}(e^{(p^{k+1})})^{p-1}\ldots (e^{(p^{k+n-1})})^{p-1}$,\
\item $(e^{(p^k)})^p=0=(f^{(p^k)})^p$,\
\item $X_k^p=X_k$
\end{enumerate}for all $k\geq0$ and $n>0$. Here $X_k:=[e^{(p^k)},f^{(p^k)}]$.\end{theorem}

\begin{remarks}\
\begin{enumerate}\item Relation 1 says that the Lie subalgebra of $H$ generated by $e^{(p^k)}$ and $f^{(p^k)}$ is isomorphic to $\mathfrak{sl}_2$. Relations 5 and 6 say that the subalgebra generated by this Lie subalgebra is in fact the restricted enveloping algebra. Relations 2, 3 and 4 indicate how these copies of the restricted enveloping algebra fit together.\
\item Notice that the Frobenius-splitting follows directly from these relations. Indeed there is a map $\theta_T:H\to H$ of algebras given by sending $e^{(p^k)}\mapsto e^{(p^{k+1})}$, $f^{(p^k)}\mapsto f^{(p^{k+1})}$ for all $k\geq 0$. This is a right inverse to $Fr$, and its image is the required subalgebra.\
\item The choice of basis $e,f,h$ depends not only on $T$ but also on a choice of Borel subgroup $B$ containing $T$. However, the map $\theta_T$ defined above does not, which justifies the notation.\end{enumerate}\end{remarks}

\noindent The proof of Theorem \ref{thm:1} is composed of two parts. First we demonstrate that, assuming relations 1,2,3,4,5 and 6, $H$ is generated by $e,e^{(p)},e^{(p^2)},\ldots$ and $f,f^{(p)},f^{(p^2)},\ldots$, subject to those relations. We then prove the relations, of which all but 6 are very easy. The proof is completely elementary.\

\section{Preliminaries}

\textbf{Kostant's $\mathbb{Z}$-form} (see \cite{J}, chapters 10, 11). We present an analogue of the PBW theorem which holds for in particular for reductive algebraic groups $G$. For simplicity we treat the case $G=SL_2$.\

Consider the algebraic group $(SL_2)_\mathbb{Z}$, flat over $\mathbb{Z}$. Its base-change to $\mathbb{F}_p$ is the algebraic group $SL_2=(SL_2)_{\mathbb{F}_p}$ over $\mathbb{F}_p$. Its base change to $\mathbb{Q}$ is the algebraic group $(SL_2)_\mathbb{Q}$ over $\mathbb{Q}$. The integral distribution algebra $Dist((SL_2)_\mathbb{Z})$ is free over $\mathbb{Z}$, and we have the identifications: \begin{align*}  \begin{matrix}
Dist((SL_2)_\mathbb{Q}) & = & Dist((SL_2)_\mathbb{Z})\otimes_\mathbb{Z}\mathbb{Q}\\ 
H & = & Dist((SL_2)_\mathbb{Z})\otimes_\mathbb{Z}\mathbb{F}_p
\end{matrix}\end{align*}This is Kostant's $\mathbb{Z}$-form.\
We have also similar compatibilities between the Lie algebras, and the chosen basis $e,f,h$ of $(\mathfrak{sl_2})_{\mathbb{F}_p}$ lifts to a standard basis, abusively also denoted $e,f,h$, of $(\mathfrak{sl_2})_{\mathbb{Z}}$. Now $Dist((SL_2)_\mathbb{Q})$ is nothing more that the universal enveloping algebra of $(\mathfrak{sl}_2)_\mathbb{Q}$. Thus in order to give a basis (together with structure constants) for $H$, it suffices to do so for $Dist((SL_2)_\mathbb{Z})$ (then reduce modulo $p$); and in order to do so for $Dist((SL_2)_\mathbb{Z})$, it suffices to present it as a certain integral form of the (rational) universal enveloping algebra. Indeed, we have:\begin{align*}Dist((SL_2)_\mathbb{Z})=span_\mathbb{Z}\left\{\frac{f^a}{a!}.\binom{h}{b}.\frac{e^c}{c!}\right\}_{a,b,c\in\mathbb{Z}_{\geq0}}\subset U((\mathfrak{sl_2})_\mathbb{Q}).\end{align*} We will write $f^{(a)}=f^a/a!, h^{(b)}=\binom{h}{b}, e^{(c)}=e^c/c!$, and denote their images in $H$ the same way. Observe that $e^p=p!e^{(p)}=0$ in $H$, and similarly $h^p=h$ and $f^p=0$ in $H$. We have the following identity (which essentially determines the structure constants) in $Dist((SL_2)_\mathbb{Z})$ (and hence in $H$):

\begin{lemma}$e^{(r)}f^{(s)}=\sum_{k=0}^{\infty}f^{(s-k)}\binom{h-s-r+2k}{k}e^{(r-k)}$.\end{lemma}

Here, by definition $f^{(a)}=0=e^{(a)}$ for any $a<0$. Also $\binom{h+a}{k}=\sum_{j=0}^{k}\binom{h}{j}\binom{a}{k-j}$ remains in Kostant's $\mathbb{Z}$-form, and thus makes sense as an element of $H$.\\

\noindent\textbf{The Casimir element} Recall that the center of $U((\mathfrak{sl}_2)_\mathbb{Q})$ is the polynomial subalgebra generated by $\delta=4fe+(h+1)^2=4ef+(h-1)^2$. For technical reasons we may prefer to replace the base $\mathbb{Z}$ in the above considerations by $\mathbb{Z}_{(p)}$ (its localization at $p$). Since $p\neq 2$, we may thus consider $\delta/4$ as an element of the `integral' form $Dist((SL_2)_{\mathbb{Z}_{(p)}})$ of $H$.\\

\noindent We are now ready to begin the proof. The real meat is in Section \ref{check}; the reader may wish to skip there directly.

\section{Sufficiency of the relations}

\begin{proposition} Assume that relations 1,2,3,4,5 and 6 hold. Then $H$ is generated by $e,e^{(p)},e^{(p^2)},\ldots$ and $f,f^{(p)},f^{(p^2)},\ldots$, subject to (only) those relations.

\begin{proof}Let $H'$ denote the algebra generated by the symbols $e,e^{(p)},e^{(p^2)},\ldots$ and $f,f^{(p)},f^{(p^2)},\ldots$, subject to relations 1,2,3,4,5 and 6. This is not intended as a subalgebra of $H$, but rather an abstract algebra. By assumption there is an obvious map $H'\to H$; we have to show that this is an isomorphism.\begin{lemma}\label{lem2} \noindent Every element of $H'$ is a linear combination of elements of the form \begin{align*}(f)^{a_0}\ldots (f^{(p^n)})^{a_n} X_0^{b_0}\ldots X_n^{b_n}(e^{(p^n)})^{c_n}\ldots (e)^{c_0}.\end{align*}\begin{proof} Since $H'$ is generated by $e,e^{(p)},e^{(p^2)},\ldots$ and $f,f^{(p)},f^{(p^2)},\ldots$, it suffices to show that every monomial in these generators may be expressed as above. Let $\xi=\xi_1\ldots \xi_t$, with each $\xi_i=e^{(p^k)}$ or $f^{(p^k)}$ for some $k$, be such a monomial. We define the \emph{weight} of such a monomial to be the sum of the formal exponents of its factors, and the \emph{disorder} of such a monomial to be the number of pairs of factors which are out of order, i.e. the number of pairs $i<j$ with $\xi_i=e^{(p^k)}$ and $\xi_j=f^{(p^l)}$ for some $k,l$. Weight and disorder are both non-negative integers.\\

\noindent If $\xi$ has zero weight or zero disorder, then it is already of the required form. So assume both these quantities are positive; we proceed by induction on (weight, disorder) in $\mathbb{Z}_{\geq 0}\times \mathbb{Z}_{\geq 0}$, ordered lexicographically. Since $\xi$ has positive disorder, the set $\{i<j:\xi_i=e^{(p^k)},\xi_j=f^{(p^l)}\operatorname{ for ~some }k,l\}$ is non-empty. Choose an element $i<j$ which minimizes $\operatorname{min}(k,l)$. Assume that $k\leq l$; the other case is similar. Then for any factor $\xi_i=e^{(p^r)}$ with $r<k$, we may use relation 3 to move $\xi_i$ to the right-hand side of $\xi$; likewise for any factor $\xi_i=f^{(p^r)}$ with $r<k$, we may use relation 3 to move $\xi_i$ to the left-hand side of $\xi$. Hence $\xi=\alpha\xi'\beta$ where $\alpha$ is a monomial in $f,\ldots,f^{(p^{k-1})}$, $\beta$ is a monomial in $e,\ldots,e^{(p^{k-1})}$ and $\xi'$ is a monomial in $e^{(p^k)},e^{(p^{k+1})},\ldots$ and $f^{(p^k)},f^{(p^{k+1})},\ldots$ with some factor $e^{(p^k)}$ appearing to the left of some factor $f^{(p^l)}$, for some $l\geq k$. $\xi'$ has lower weight than $\xi$, and hence is less than $\xi$ in the lexicographic order, unless $\xi'=\xi$. So assume $\xi'=\xi$ (else done).\\

\noindent Recall we have $i<j$ with $\xi_i=e^{(p^k)}$ and $\xi_j=f^{(p^l)}$. Let $m>i$ be minimal such that $\xi_m=f^{(p^r)}$ for some $r$. Then we can reorder the factors of $\xi$ so that $\xi_{m-1}=e^{(p^k)}$. In other words, we reduce to the case $\xi=\xi_1\ldots\xi_{m-2}e^{(p^k)}f^{(p^l)}\xi_{m+1}\ldots \xi_t$, with each $\xi_i$ being one of $e^{(p^k)},e^{(p^{k+1})},\ldots$ or $f^{(p^k)},f^{(p^{k+1})},\ldots$, and $l\geq k$. Then: \begin{align*}\xi  = \xi_1\ldots\xi_{m-2}f^{(p^l)}e^{(p^k)}\xi_{m+1}\ldots \xi_t + \xi_1\ldots\xi_{m-2}[e^{(p^k)},f^{(p^l)}]\xi_{m+1}\ldots \xi_t.\end{align*} The first summand of the RHS has the same weight as $\xi$, but lower disorder, so is less than $\xi$ in the lexicographic ordering and may be ignored (by induction). The second summand is calculated using relation 1 or relation 4, depending on the value of $l$. If $l=k$, then it is equal to $\xi_1\ldots\xi_{m-2}X_k\xi_{m+1}\ldots \xi_t$. By relations 1 and 2, this is equal to $\xi_1\ldots\xi_{m-2}\xi_{m+1}\ldots \xi_t(X_k+2q)$, where $q$ is the difference between the number of factors equal to $e^{(p^k)}$, and the number of factors equal to $f^{(p^k)}$, amongst $\xi_{m+1},\ldots, \xi_t$. Otherwise, $l=k+r>k$ and the second summand is equal to $(-1)^r\xi_1\ldots\xi_{m-2}(f^{(p^k)})^{p-1}(f^{(p^{k+1})})^{p-1}\ldots (f^{(p^{k+r-1})})^{p-1}\xi_{m+1}\ldots \xi_t(X_k+1+2q)$. In either case, we see that the second summand is equal to $\pm\xi'(X_k+c)$ for some monomial $\xi'$ in the elements $e^{(p^k)},e^{(p^{k+1})},\ldots$, $f^{(p^k)},f^{(p^{k+1})},\ldots$, of lower weight than $\xi$, and some constant $c$.\\

\noindent Note that the subalgebra of $H'$ generated by $e^{(p^k)},e^{(p^{k+1})},\ldots$, $f^{(p^k)},f^{(p^{k+1})},\ldots$ is isomorphic to $H'$, via $e^{(p^r)}\mapsto e^{(p^{r+k})}$, $f^{(p^r)}\mapsto f^{(p^{r+k})}$. Let $\xi''$ be the preimage of $\xi'$ under this map; its weight is at most that of $\xi'$, and so by induction we may write it as a linear combination of elements of the form \begin{align*}(f)^{a_0}\ldots (f^{(p^n)})^{a_n} X_0^{b_0}\ldots X_n^{b_n}(e^{(p^n)})^{c_n}\ldots (e)^{c_0}.\end{align*}  Thus $\xi'$ is written as a linear combination of elements of the form \begin{align*}(f^{(p^k)})^{a_0}\ldots (f^{(p^{k+n})})^{a_n} X_k^{b_0}\ldots X_{k+n}^{b_n}(e^{(p^{k+n})})^{c_n}\ldots (e^{(p^k)})^{c_0}.\end{align*} We conclude by observing that \begin{align*}& (f^{(p^k)})^{a_0}\ldots (f^{(p^{k+n})})^{a_n} X_k^{b_0}\ldots X_{k+n}^{b_n}(e^{(p^{k+n})})^{c_n}\ldots (e^{(p^k)})^{c_0}(X_k+c)=\\
&  (f^{(p^k)})^{a_0}\ldots (f^{(p^{k+n})})^{a_n} X_k^{b_0+1}\ldots X_{k+n}^{b_n}(e^{(p^{k+n})})^{c_n}\ldots (e^{(p^k)})^{c_0}\\
& +(c-2c_0)(f^{(p^k)})^{a_0}\ldots (f^{(p^{k+n})})^{a_n} X_k^{b_0}\ldots X_{k+n}^{b_n}(e^{(p^{k+n})})^{c_n}\ldots (e^{(p^k)})^{c_0}\end{align*} has the required form (note that relation 2 implies that all $X_i$ commute).\end{proof}\end{lemma}

\noindent Note that relations 5 and 6 allow us to take the exponents $a_i,b_i,c_i<p$ in the statement of Lemma \ref{lem2}.\\

\noindent Now we show that the elements $(f)^{a_0}\ldots (f^{(p^n)})^{a_n} X_0^{b_0}\ldots X_n^{b_n}(e^{(p^n)})^{c_n}\ldots (e)^{c_0}$ with $a_i,b_i,c_i<p$ of $H$ form a basis; then we will be done. We know that $H$ has a basis consisting of \begin{align*}(f)^{a_0}\ldots (f^{(p^n)})^{a_n} \binom{h}{p^0}^{b_0}\ldots \binom{h}{p^n}^{b_n}(e^{(p^n)})^{c_n}\ldots (e)^{c_0}\end{align*} with $a_i,b_i,c_i<p$. Define the \emph{weight} of such a monomial to be $b_0+b_1p+\ldots+b_np^n$. This induces a partial order on the given basis, where one monomial is less than another if it has lower weight. In fact, the given basis is well (partially) ordered for this ordering. Then it is easy to see that \begin{align*}(f)^{a_0}\ldots (f^{(p^n)})^{a_n} X_0^{b_0}\ldots X_n^{b_n}(e^{(p^n)})^{c_n}\ldots (e)^{c_0} & =  (f)^{a_0}\ldots (f^{(p^n)})^{a_n} \binom{h}{p^0}^{b_0}\ldots \binom{h}{p^n}^{b_n}(e^{(p^n)})^{c_n}\ldots (e)^{c_0}\\
& + \operatorname{linear~ combination~ of~ lower~ weight~ monomials}.\end{align*}
from which it follows that the linear map $H\to H$ given by mapping \begin{align*}(f)^{a_0}\ldots (f^{(p^n)})^{a_n} \binom{h}{p^0}^{b_0}\ldots \binom{h}{p^n}^{b_n}(e^{(p^n)})^{c_n}\ldots (e)^{c_0}\mapsto (f)^{a_0}\ldots (f^{(p^n)})^{a_n} X_0^{b_0}\ldots X_n^{b_n}(e^{(p^n)})^{c_n}\ldots (e)^{c_0}\end{align*} is an isomorphism.\end{proof}\end{proposition}

\section{Checking the relations}\label{check}

\noindent It remains to prove that relations 1,2,3,4,5 and 6 hold in $H$. Relations 3 and 5 are trivial. Relations 1,2 and 4 are short calculations:

\begin{lemma}Relation 1 holds in $H$.
\begin{proof}We have $X_k=[e^{(p^k)},f^{(p^k)}]=\sum_{i=1}^{p^k}\binom{h}{i}f^{(p^k-i)}e^{(p^k-i)}$, so that \begin{align*}[X_k,e^{(p^k)}] & = \sum_{i=1}^{p^k}\left[\binom{h}{i}f^{(p^k-i)}e^{(p^k-i)},e^{(p^k)}\right]\\
& = -\sum_{i=1}^{p^k-1}\binom{h}{i}[e^{(p^k)},f^{(p^k-i)}]e^{(p^k-i)}+\left[\binom{h}{p^k},e^{(p^k)}\right]\\
& = -\sum_{i=1}^{p^k-1}\binom{h}{i}\sum_{j=1}^{p^k-i}f^{(p^k-i-j)}\binom{h+i+2j}{j}e^{(p^k-j)}e^{(p^k-i)}+\left[\binom{h}{p^k},e^{(p^k)}\right]\\
& = \left[\binom{h}{p^k},e^{(p^k)}\right]\\
& = \left\{\binom{h}{p^k}-\binom{h-2p^k}{p^k}\right\}e^{(p^k)}\\
& = 2e^{(p^k)}\end{align*}as required. Similarly $[X_k,f^{(p^k)}]=-2f^{(p^k)}$.\end{proof}\end{lemma}

\begin{lemma}Relation 2 holds in $H$.
\begin{proof}We have \begin{align*}[X_k,e^{(p^{k+r})}] & = \sum_{i=1}^{p^k}\left[\binom{h}{i}f^{(p^k-i)}e^{(p^k-i)},e^{(p^{k+r})}\right]\\
& = -\sum_{i=1}^{p^k}\binom{h}{i}[e^{(p^{k+r})},f^{(p^k-i)}]e^{(p^k-i)}\\
& = -\sum_{i=1}^{p^k-1}\binom{h}{i}\sum_{j=1}^{p^k-i}f^{(p^k-i-j)}\binom{h+i+2j}{j}e^{(p^{k+r}-j)}e^{(p^k-i)}\\
& = 0
\end{align*}as required. Similarly $[X_k,f^{(p^{k+r})}]=0$.\end{proof}\end{lemma}

\begin{lemma}Relation 4 holds in $H$.
\begin{proof}We have \begin{align*}[e^{(p^{k+r})},f^{(p^{k})}] & = \sum_{i=1}^{p^{(k)}}f^{(p^k-i)}\binom{h-p^k+2i}{i}e^{(p^{k+r}-i)}\\
& = \sum_{i=1}^{p^{(k)}}\binom{h+p^k}{i}f^{(p^k-i)}e^{(p^{k+r}-i)}\\
& = \sum_{i=1}^{p^{(k)}}\binom{h}{i}f^{(p^k-i)}e^{(p^{k+r}-i)}+e^{(p^{k+r}-p^k)}\\
& = \sum_{i=1}^{p^{(k)}}\binom{h}{i}f^{(p^k-i)}e^{(p^{k}-i)}e^{(p^{k+r}-p^k)}\binom{p^{k+r}-i}{p^k-i}^{-1}+e^{(p^{k+r}-p^k)}\\
& = \sum_{i=1}^{p^{(k)}}\binom{h}{i}f^{(p^k-i)}e^{(p^{k}-i)}e^{(p^{k+r}-p^k)}+e^{(p^{k+r}-p^k)}\\
& = (X_k+1)e^{(p^{k+r}-p^k)}\\
& = (-1)^{r}(X_k+1)(e^{(p^k)})^{p-1}(e^{(p^{k+1})})^{p-1}\ldots (e^{(p^{k+r-1})})^{p-1}\end{align*}
as required. Similarly $[e^{(p^k)},f^{(p^{k+r})}]=(-1)^{r}(f^{(p^k)})^{p-1}(f^{(p^{k+1})})^{p-1}\ldots (f^{(p^{k+r-1})})^{p-1}(X_k+1)$.\end{proof}\end{lemma}

\noindent Now set $t_k:=X_k-\binom{h}{p^k}\in H$. Then relation 6 is equivalent to the statement that $t_k^p=t_k$. In fact, we prove the following

\begin{theorem}$t_k^2=t_k$.
\begin{proof}We first prove the case $k=1$ (case $k=0$ is trivial). To that end, let $H'$ denote the Kostant $\mathbb{Z}$-form and $H''$ denote $H'\otimes_{\mathbb{Z}}\mathbb{Z}_{(p)}$, so that $H=H''\otimes_{\mathbb{Z}_{(p)}}\mathbb{F}_p$. We will construct a lift of $X_1$ to $H'$. Denote the central (Casimir) element $4fe+(h+1)^2=4ef+(h-1)^2\in H''$ by $\delta$. Then in $H''$ we have the following equalities:\begin{align*} 4^p(p-1)!^2e^{(p)}f^{(p)}=(\delta-(h-1)^2)(\delta-(h-3)^2)\ldots(\delta-(h-2p+3)^2)(\delta-(h-2p+1)^2)/p^2\\
4^p(p-1)!^2f^{(p)}e^{(p)}=(\delta-(h+1)^2)(\delta-(h+3)^2)\ldots(\delta-(h+2p-3)^2)(\delta-(h+2p-1)^2)/p^2.\end{align*}The difference between the above expressions is a degree $p-1$ polynomial in $\delta$ with coefficients in $\frac{1}{p^2}\mathbb{Z}[h]$; call it $Q=Q_{p-1}\delta^{p-1}+Q_{p-2}\delta^{p-2}+\ldots+Q_1\delta+Q_0\in H''$. Notice that for any $m$, $\delta^m=4^mf^me^m+\chi_{k-1}f^{m-1}e^{m-1}+\ldots+\chi_0$ for some $\chi_i\in\mathbb{Z}[h]$. Now is the crux: it follows (by descending induction on $i$) that, for each $i$, $Q_i\in H''$. So the image of $Q$ in $H$ is equal to $\overline{Q}=\overline{Q_{p-1}}\delta^{p-1}+\overline{Q_{p-2}}\delta^{p-2}+\ldots+\overline{Q_1}\delta+\overline{Q_0}$. Here $\overline{Q_i}$ stands for the image of $Q_i$ in $H$; it is an element of the distribution algebra of maximal torus $T$. By an abuse of notation $\delta$ stands for its image in $H$. Of course, $\overline{Q}=4X_1$.\\

\noindent Observe that $\delta^p-2\delta^{\frac{p+1}{2}}+\delta=\prod_{j\in\mathbb{F}_p}(\delta-j^2)=0$ in $H$; this is the minimal polynomial of $\delta$. Likewise $h^p-h$ is the minimal polynomial of $h$ in $H$. Thus in particular the subalgebra of $H$ generated by $h,\delta$ is isomorphic to \begin{align*}\mathbb{F}_p[h]/(h^p-h)\otimes\mathbb{F}_p[\delta]/(\delta^p-2\delta^{\frac{p+1}{2}}+\delta)\cong \left\{\prod_{i\in\mathbb{F}_p}\mathbb{F}_p\right\}\otimes\left\{\prod_{j^2=0}\mathbb{F}_p\times \prod_{j^2\in\mathbb{F}_p^\times}\mathbb{F}_p[\epsilon]/(\epsilon^2)\right\}.\end{align*}

\noindent Here the map from $\mathbb{F}_p[h,\delta]$ to the $i,j^2$ factor $\mathbb{F}_p[\epsilon]/(\epsilon^2)$ sends $h$ to $i$ and $\delta$ to $j^2+\epsilon$ (for $j^2\in\mathbb{F}_p^\times$), while the map to the $i,0$ factor sends $h$ to $i$ and $\delta$ to $0$.\\

\noindent We know that $t_1=X_1-\binom{h}{p}=\sum_{k=1}^{p-1}f^{(p-k)}e^{(p-k)}\binom{h}{k}\in\mathbb{F}_p[h,\delta]\cong\mathbb{F}_p[h]\otimes\mathbb{F}_p[\delta]$, from which it follows that $\overline{Q_1},\ldots,\overline{Q}_{p-1}\in\mathbb{F}_p[h]$, while $\overline{Q_0}-4\binom{h}{p}\in\mathbb{F}_p[h]$. Thus we have $4t_1=\overline{Q_{p-1}}\delta^{p-1}+\overline{Q_{p-2}}\delta^{p-2}+\ldots+\overline{Q_1}\delta+\overline{Q_0}-4\binom{h}{p}$, and to check that $t_1^2=t_1$, it suffices to check that, for each $i,j^2$, the image of $4t_1$ in the $i,j^2$ factor above is equal to $0$ or $4$.\\

\noindent First we should check that for $j^2\in\mathbb{F}_p^\times$, and any $i$, the image of $4t_1$ in the $i,j^2$ factor is constant (its coefficient of $\epsilon$ is $0$). So assume $j^2$ is a non-zero quadratic residue in $\mathbb{F}_p$. Choose any lift $\tilde{j}$ of $j$ to $\mathbb{Z}$. Write \begin{align*}{Q}={Q_{p-1}}\delta^{p-1}+{Q_{p-2}}\delta^{p-2}+\ldots+{Q_1}\delta+{Q_0}={R_{p-1}}(\delta-\tilde{j}^2)^{p-1}+{R_{p-2}}(\delta-\tilde{j}^2)^{p-2}+\ldots+{R_1}(\delta-\tilde{j}^2)+{R_0}\end{align*} for some $R_i\in\frac{1}{p^2}\mathbb{Z}[h]\cap H''$. We need to show that $\overline{R_1}=0$. It is equivalent to showing that $R_1/p\in  H''$, or equivalently that $p^2R_1\in \mathbb{Z}[h]$ maps every integer value of $h$ to an element of $p^3\mathbb{Z}_{(p)}$. \\


\noindent So fix any value of $h\in\mathbb{Z}$. $p^2R_1\in \mathbb{Z}[h]$ is the coefficient of $\delta-\tilde{j}^2$ in the $\delta-\tilde{j}^2$-adic expansion of \begin{align*} & (\delta-(h-1)^2)(\delta-(h-3)^2)\ldots(\delta-(h-2p+3)^2)(\delta-(h-2p+1)^2)\\
- & (\delta-(h+1)^2)(\delta-(h+3)^2)\ldots(\delta-(h+2p-3)^2)(\delta-(h+2p-1)^2).\end{align*} So it is the difference between the coefficients of $\delta-\tilde{j}^2$ in the $\delta-\tilde{j}^2$-adic expansions of \begin{align*} & (\delta-(h-1)^2)(\delta-(h-3)^2)\ldots(\delta-(h-2p+3)^2)(\delta-(h-2p+1)^2)\\
= & \prod_{l=1}^{p}(\delta-\tilde{j}^2+(\tilde{j}+h-2l+1)(\tilde{j}-h+2l-1))\end{align*} and \begin{align*} & (\delta-(h+1)^2)(\delta-(h+3)^2)\ldots(\delta-(h+2p-3)^2)(\delta-(h+2p-1)^2)\\
= & \prod_{l=1}^{p}(\delta-\tilde{j}^2+(\tilde{j}+h+2l-1)(\tilde{j}-h-2l+1)).\end{align*} Let us denote the former coefficient by $\chi(h)$; then the latter coefficient is equal to $\chi(h+2p)$. We have \begin{align*}\chi(h) & =\sum_{i=1}^{p}\prod_{\substack{1\leq l\leq p\\
l\neq i}}(\tilde{j}+h-2l+1)(\tilde{j}-h+2l-1)\\
& = \frac{1}{2\tilde{j}}\sum_{i=1}^{p}(\tilde{j}+h-2i+1+\tilde{j}-h+2i-1)\prod_{\substack{1\leq l\leq p\\
l\neq i}}(\tilde{j}+h-2l+1)(\tilde{j}-h+2l-1)\\
& = \frac{1}{2\tilde{j}}\sum_{i=1}^{p}(\tilde{j}+h-2i+1+\tilde{j}-h+2i-1)\prod_{\substack{1\leq l\leq p\\
l\neq i}}(\tilde{j}+h-2l+1)(\tilde{j}-h+2l-1)\\
& = \frac{1}{2\tilde{j}}\left\{\sum_{i=1}^{p}\prod_{\substack{1\leq l\leq p\\
l\neq i}}(\tilde{j}+h-2l+1).\prod_{1\leq l\leq p}(\tilde{j}-h+2l-1)+\sum_{i=1}^{p}\prod_{1\leq l\leq p}(\tilde{j}+h-2l+1).\prod_{\substack{1\leq l\leq p\\
l\neq i}}(\tilde{j}-h+2l-1)\right\}\end{align*}For each $1\leq i\leq p$, there exists a unique $1\leq\tau(i)\leq p$ such that $\tilde{j}-i+\tau(i)\equiv 0 \operatorname{mod~}p$; $\tau$ is a bijection. Note that $(\tilde{j}+h-2i+1)+(\tilde{j}-h+2\tau(i)-1)=2(\tilde{j}-i+\tau(i))$. So we have \begin{align*}\chi(h) & = \frac{1}{2\tilde{j}}\sum_{i=1}^{p}2(\tilde{j}-i+\tau(i))\prod_{\substack{1\leq l\leq p\\
l\neq i}}(\tilde{j}+h-2l+1).\prod_{\substack{1\leq l\leq p\\
l\neq\tau(i)}}(\tilde{j}-h+2l-1)\\
& = \frac{1}{\tilde{j}}\sum_{i=1}^{p}(\tilde{j}-i+\tau(i))\prod_{\substack{1\leq l\leq p\\
l\neq i}}(\tilde{j}+h-2l+1)(\tilde{j}-h+2\tau(l)-1)\end{align*}There is a unique $l_0$, $1\leq l_0\leq p$, such that $\tilde{j}+h-2l_0+1\equiv0\operatorname{mod~}p$. Then $\tau(l_0)$ is the unique integer between $1$ and $p$ such that $\tilde{j}-h+2\tau(l_0)-1\equiv0\operatorname{mod~}p$. Since $\tilde{j}$ is not divisible by $p$, it follows that for every $i\neq l_0$ with $1\leq i\leq p$, the corresponding summand above is divisible by $p^3$. So set \begin{align*}\phi(h)=\prod_{\substack{1\leq l\leq p\\
l\neq l_0}}(\tilde{j}+h-2l+1)(\tilde{j}-h+2\tau(l)-1);\end{align*}we need to show that $\phi(h)-\phi(h+2p)$ is divisible by $p^2$, or equivalently, that $\phi'(h)$ is divisible by $p$. But we have \begin{align*}\phi'(h)=\sum_{\substack{1\leq i\leq p\\
i\neq l_0}}\prod_{\substack{1\leq l\leq p\\
l\neq i,l_0}}(\tilde{j}+h-2l+1).\prod_{\substack{1\leq l\leq p\\
l\neq l_0}}(\tilde{j}-h+2\tau(l)-1)-\sum_{\substack{1\leq i\leq p\\
i\neq l_0}}\prod_{\substack{1\leq l\leq p\\
l\neq l_0}}(\tilde{j}+h-2l+1).\prod_{\substack{1\leq l\leq p\\
l\neq i,l_0}}(\tilde{j}-h+2\tau(l)-1)\end{align*}As $i$ ranges from $1$ to $p$, excluding $l_0$, the expressions $\tilde{j}+h-2l+1$, $\tilde{j}-h+2\tau(l)-1$ both take each non-zero residue modulo $p$ precisely once. Therefore \begin{align*}\phi'(h)\equiv \sum_{i=1}^{p-1}\frac{(p-1)!}{i}(p-1)!-\sum_{i=1}^{p-1}(p-1)!\frac{(p-1)!}{i}\equiv 0\operatorname{mod~}p,\end{align*}as required.\\

\noindent Now we need to check that for any $i,j^2$, the image of $4t_1$ in $\mathbb{F}_p[h,\delta]/(h-i,\delta-j^2)$ is $0$ or $4$. This is proved similarly. Indeed, choose any lift $\tilde{j}$ of $j$, and let $\tilde{i}$ be the unique lift of $i$ such that $0\leq\tilde{i}<p$ (so that $\binom{i}{p}=0$); we should check that $Q_{p-1}(\tilde{i})\tilde{j}^{p-1}+Q_{p-2}(\tilde{i})\tilde{j}^{p-2}+\ldots+Q_1(\tilde{i})\tilde{j}+Q_0(\tilde{i})$, which is an integer, is congruent to $0$ or $4$ modulo $p$. Equivalently we should show that \begin{align*} & (\tilde{j}^2-(\tilde{i}-1)^2)(\tilde{j}^2-(\tilde{i}-3)^2)\ldots(\tilde{j}^2-(\tilde{i}-2p+3)^2)(\tilde{j}^2-(\tilde{i}-2p+1)^2)/p^2\\
- & (\tilde{j}^2-(\tilde{i}+1)^2)(\tilde{j}^2-(\tilde{i}+3)^2)\ldots(\tilde{j}^2-(\tilde{i}+2p-3)^2)(\tilde{j}^2-(\tilde{i}+2p-1)^2)/p^2\end{align*} is congruent to $0$ or $4$ modulo $p$. Let $a$, $b$ be the unique integers between $1$ and $p$ such that $\tilde{j}-\tilde{i}+2a-1$, $\tilde{j}+\tilde{i}-2b+1$ are both divisible by $p$, and write them respectively as $rp$, $sp$. Then we need only show that $rs-(r-2)(s+2)=-2(r-s)+4$ is congruent to $0$ or $4$ modulo $p$, or equivalently that $r-s$ is congruent to $0$ or $2$ modulo $p$. But $(r-s)p=(2a-1)+(2b-1)-2\tilde{i}$ is an even multiple of $p$ satisfying $-2p+4\leq (r-s)p \leq 4p-2$, so is equal to $0$ or $2p$.\\

\noindent This proves that $t_1^2=t_1$. We show inductively that $t_k^2=t_k$. We have \begin{align*}X_k & =\sum_{i=1}^{p^k}\binom{h}{i}f^{(p^k-i)}e^{(p^k-i)}\\
& = \sum_{j=1}^{p}\sum_{i=1}^{p^{k-1}}\binom{h}{p^{k-1}(j-1)+i}f^{(p^{k-1}(p-j)+p^{k-1}-i)}e^{(p^{k-1}(p-j)+p^{k-1}-i)}\\
& = \sum_{j=1}^{p}\sum_{i=1}^{p^{k-1}-1}\binom{h}{p^{k-1}(j-1)}\binom{h}{i}f^{(p^{k-1}(p-j))}f^{(p^{k-1}-i)}e^{(p^{k-1}-i)}e^{(p^{k-1}(p-j))}+\sum_{j=1}^{p}\binom{h}{p^{k-1}j}f^{(p^{k-1}(p-j))}e^{(p^{k-1}(p-j))}\\
& = \left(X_{k-1}-\binom{h}{p^{k-1}}\right)\sum_{j=1}^{p}\binom{h}{p^{k-1}(j-1)}f^{(p^{k-1}(p-j))}e^{(p^{k-1}(p-j))}+\sum_{j=1}^{p}\binom{h}{p^{k-1}j}f^{(p^{k-1}(p-j))}e^{(p^{k-1}(p-j))}\\
& = t_{k-1}\sum_{j=1}^{p}\binom{X_{k-1}-t_{k-1}}{j-1}f^{(p^{k-1}(p-j))}e^{(p^{k-1}(p-j))}+\sum_{j=1}^{p-1}\binom{X_{k-1}-t_{k-1}}{j}f^{(p^{k-1}(p-j))}e^{(p^{k-1}(p-j))}+\binom{h}{p^k}\end{align*}so that \begin{align*}t_k & = t_{k-1}\binom{X_{k-1}-t_{k-1}}{p-1} + \sum_{j=1}^{p-1}\left(t_{k-1}\binom{X_{k-1}-t_{k-1}}{j-1}+\binom{X_{k-1}-t_{k-1}}{j}\right)f^{(p^{k-1}(p-j))}e^{(p^{k-1}(p-j))}\\
& = t_{k-1}\binom{X_{k-1}-t_{k-1}}{p-1} + \sum_{j=1}^{p-1}\binom{X_{k-1}}{j}f^{(p^{k-1}(p-j))}e^{(p^{k-1}(p-j))}\end{align*}since $t_{k-1}^2=t_{k-1}$. Moreover since $X_{k-1}^p=X_{k-1}$, it follows that the subalgebra generated by $e^{(p^{k-1})}$, $f^{(p^{k-1})}$ is the restricted enveloping algebra. We have already proved that \begin{align*}\sum_{j=1}^{p-1}\binom{h}{j}f^{(p-j)}e^{(p-j)}=\sum_{j=1}^{p-1}\binom{h}{j}f^{p-j}e^{p-j}/(p-j)!^2\end{align*} is idempotent, since it is equal to $t_1$. Thus also \begin{align*}\sum_{j=1}^{p-1}\binom{X_{k-1}}{j}f^{(p^{k-1}(p-j))}e^{(p^{k-1}(p-j))}=\sum_{j=1}^{p-1}\binom{X_{k-1}}{j}(f^{(p^{k-1})})^{p-j}(e^{(p^{k-1})})^{p-j}/(p-j)!^2\end{align*} is idempotent. Since $X_{k-1}-t_{k-1}$ is fixed under raising to the $p^{th}$ power, we have $(X_{k-1}-t_{k-1})\binom{X_{k-1}-t_{k-1}}{p-1}=-\binom{X_{k-1}-t_{k-1}}{p-1}$, so $\binom{X_{k-1}-t_{k-1}}{p-1}^2=\binom{-1}{p-1}\binom{X_{k-1}-t_{k-1}}{p-1}=\binom{X_{k-1}-t_{k-1}}{p-1}$ is idempotent. Therefore $t_{k-1}\binom{X_{k-1}-t_{k-1}}{p-1}$ is idempotent since $t_{k-1}$ commutes with $X_{k-1}$. Finally, $t_{k-1}\binom{X_{k-1}-t_{k-1}}{p-1}X_k=t_{k-1}\binom{X_{k-1}-t_{k-1}}{p-1}(t_k-1)=(t_{k-1}^2-t_{k-1})\binom{X_{k-1}-t_{k-1}}{p-1}=0$, and $\sum_{j=1}^{p-1}\binom{X_{k-1}}{j}f^{(p^{k-1}(p-j))}e^{(p^{k-1}(p-j))}$ is divisible (on the left) $X_{k-1}$ so that the idempotents $t_{k-1}\binom{X_{k-1}-t_{k-1}}{p-1}$ and $\sum_{j=1}^{p-1}\binom{X_{k-1}}{j}f^{(p^{k-1}(p-j))}e^{(p^{k-1}(p-j))}$ are orthogonal.\end{proof}\end{theorem}






\begin{thebibliography}{100}

\bibitem{J} Jens Carsten Jantzen, \textit{Representations of Algebraic Groups}, Second Edition, Amer. Math, Soc., Providence, RI, 2003, reprinted in 2007.

\bibitem{G1} Michel Gros, \textit{A Splitting of the Frobenius Morphism on the Whole Algebra of Distributions of $SL_2$}, Algebr. Represent. Theor., 2012.

\bibitem{G2} Michel Gros, Masaharu Kaneda, \textit{Contraction par Frobenius de G-modules}, arXiv:1004.1939 [math.RT], 2010.





























\end{thebibliography}
\end{document}